\documentclass{article}[11pt]
\usepackage{epsfig}
\usepackage{enumerate,graphicx}
\usepackage{amsfonts}
\usepackage{graphicx}
\usepackage{amsmath,amsthm,amssymb}
\usepackage{amssymb}
\usepackage{appendix}
\DeclareGraphicsExtensions{.pdf, .jpg, .png , .bmp, .ps}
\usepackage{tikz}
\providecommand{\U}[1]{\protect\rule{.1in}{.1in}}
\newtheorem{thm}{Theorem}[section]

\newtheorem{prop}{Proposition}[section]
\newtheorem{cor}{Corollary}[section]

\theoremstyle{remark}
\newtheorem{rem}[thm]{\bf Remark}

\newcommand{\tq}{\,:\,}\newcommand\E{{\mathbb {E}}}

\newcommand\I{{\bf 1}}

\def\1{{{\mbox{${\rm{1\negthinspace\negthinspace I}}$}}}}

\newcommand\beq{\begin{equation}}
\newcommand\eeq{\end{equation}}

\begin{document}

\title{Almost sure invariance principle for the Kantorovich distance between
the empirical and the marginal distributions  of  strong mixing sequences}

\author{J\'er\^ome Dedecker\footnote{J\'er\^ome Dedecker, Universit\'e de Paris, CNRS, MAP5, UMR 8145,
45 rue des  Saints-P\`eres,
F-75006 Paris, France.},
Florence Merlev\`ede \footnote{Florence Merlev\`ede, 
 Universit\'e Gustave Eiffel,   LAMA, UMR 8050 CNRS,  \  F-77454 Marne-La-Vall\'ee, France.}
}

\maketitle

\abstract{We  prove a strong invariance principle for the  Kantorovich distance between the empirical distribution and the
marginal distribution of  stationary  $\alpha$-mixing sequences. 
}

\medskip

\noindent{\bf Running head.} ASIP for the empirical  $W_1$
distance.

\medskip

\noindent {\bf Keywords.} Empirical process, Wasserstein distance, Almost sure invariance principle, Compact law of the iterated logarithm, Bounded law of the iterated logarithm, Conditional Value at Risk

\medskip

\noindent {\bf Mathematics Subject Classification (2010).} 60F15, 60G10, 60B12

\section{Introduction and notations}
Let $(X_i)_{i \in {\mathbb Z}}$ be a strictly stationary sequence of real-valued random variables. Define the two $\sigma$-algebras ${\mathcal F}_0= \sigma(X_i, i \leq 0)$ and ${\mathcal G}_k= \sigma(X_i, i \geq k)$, and recall that the strong mixing coefficients $(\alpha(k))_{k \geq 0}$ of Rosenblatt \cite{R} are defined by 
\begin{equation}\label{defalpha}
  \alpha(k)= \sup_{A \in {\mathcal F}_0, B \in {\mathcal G}_k} \left | {\mathbb P}(A \cap B)- {\mathbb P}(A) {\mathbb P}(B) \right | \, .
\end{equation}

Let  
$\mu$ be the common distribution of the $X_i$'s, and let
$$
   \mu_n = \frac 1 n \sum_{k=1}^n \delta_{X_k} 
$$
be the empirical measure based on $X_1, \ldots, X_n$. In this paper, we prove a strong invariance principle for the Kantorovich distance  $W_1(\mu_n, \mu)$ between $\mu_n$ and $\mu$ under a condition on the mixing coefficients $\alpha(k)$. Recall that the Kantorovich distance  (also called Wasserstein distance of order 1) between two probability measures $\mu$ and $\nu$  is defined by
$$
  W_1(\mu, \nu)= \inf_{ \pi \in M(\mu, \nu)} \int |x-y| \pi(dx, dy)  \, ,
$$
where $M(\mu, \nu)$ is the set of probability measures on ${\mathbb R}^2$ with marginals $\mu$ and $\nu$. We shall use the following well known representation for probabilities on the real line:
\begin{equation}\label{Frep}
W_1(\mu, \nu)= \int |F_\mu(x)-F_\nu(x)| dx  \, ,
\end{equation}
where $F_\mu$ is the cumulative distribution function of $\mu$. 

Let $H: t \rightarrow {\mathbb P}([X_0|>t)$ be the tail function of $|X_0|$. In the case where $(X_i)_{i \in {\mathbb Z}}$ is a sequence of independent and identically distributed (i.i.d.) random variables, del Barrio et al. \cite{dBGM} used the representation \eqref{Frep} and a general result of  Jain \cite{J} for Banach-valued random variables to prove a central limit theorem for $\sqrt n W_1(\mu_n, \mu)$.
More precisely, they showed  that $\sqrt n W_1(\mu_n, \mu)$ converges in distribution to the ${\mathbb L}_1(dt)$ norm of an ${\mathbb L}_1(dt)$-valued Gaussian random variable, provided that
\begin{equation}\label{Gine}
\int_0^\infty \sqrt{ H(t)} \ dt < \infty \, .
\end{equation} 
They also proved that $\sqrt n W_1(\mu_n, \mu)$ is stochastically bounded iff \eqref{Gine} holds, proving that this condition is necessary and sufficient for the weak convergence of $\sqrt n W_1(\mu_n, \mu)$.

Still in the i.i.d. case, we easily deduce from Chapters 8 and 10 in Ledoux and Talagrand \cite{LT} that: if \eqref{Gine} holds, then the sequence 
\begin{equation}\label{2LL}
    \frac{\sqrt{n}}{\sqrt{2 \log \log n}}W_1(\mu_n, \mu)
\end{equation}
satisfies a compact law of the iterated logarithm. 

For strongly mixing sequences in the sense of Rosenblatt \cite{R}, we proved in \cite{DM17}  the central limit theorem for  $\sqrt n W_1(\mu_n, \mu)$ under the condition 
\begin{equation} \label{condalpha}
  \int_0^{\infty} \sqrt{ \sum_{k=0}^\infty  \left ( \alpha(k) \wedge H(t)
\right ) } \ dt < \infty 
\end{equation}
 (where $a\wedge b$ means the minimum between two reals $a$ and $b$),  and we give sufficient conditions for \eqref{condalpha} to hold. Note that, in \cite{DM17}, we used a weaker version of the $\alpha$-mixing coefficients, that enables to deal with a large class of non-mixing processes in the sense of Rosenblatt \cite{R}. 

In Section \ref{MainSec}  of this paper, we prove a strong invariance principle for $W_1(\mu_n, \mu)$ under the condition \eqref{condalpha}. The compact law of the iterated logarithm for \eqref{2LL} easily follows from this strong invariance principle. In Section \ref{ES}, we  apply our general result to derive the almost sure rate of convergence of the empirical estimator of the Conditional Value at Risk ($CVaR$) for stationary $\alpha$-mixing sequences. 

In the rest of the paper, we shall use the following notation:  for two sequences $(a_n)_{n \geq 1}$ and $(b_n)_{n \geq 1}$ of positive reals, $a_n \ll b_n$ means there exists a positive constant $C$ not depending on $n$ such that $a_n \leq C b_n$ for any $n\geq 1$. 

\section{Main result} \label{MainSec}

Our main result is the following strong invariance principle for $W_1(\mu_n, \mu)$. 

\setcounter{equation}{0}

\begin{thm}\label{th1}
Assume that \eqref{condalpha} is satisfied.
Then, enlarging the probability space if necessary, there exists a sequence of i.i.d.  ${\mathbb L}_1(dt)$-valued centered  Gaussian random variables $(Z_i)_{i \geq 1}$ with covariance function defined as follows:
for any $f, g \in {\mathbb L}_\infty (dt)$,
\begin{equation}\label{covalpha}
\Gamma(f,g)= \mathrm{Cov} \left(  \int f(t) Z_1(t) \ dt , \int g(t) Z_1(t) \ dt \right )
=  \sum_{k \in {\mathbb Z}}
\iint f(t)g(s) \mathrm{Cov} (  {\bf 1}_{X_0 \leq t} ,  {\bf 1}_{X_k \leq s} )  \ ds \, dt  \, ,
\end{equation}
and such that
$$
 n W_1(\mu_n, \mu) - \int \left |  \sum_{k=1}^n Z_k(t) \right | dt = o(\sqrt{n \log \log n}) \quad
 \text{almost surely.}
$$
\end{thm}

\begin{rem}
 In \cite{Cuny17}, Cuny proved a strong invariance principle for $W_1(\mu_n, \mu)$. under the condition
\begin{equation}\label{CunyCond}
  \sum_{k=0}^\infty  \frac 1 {\sqrt {k+1}} \int_0^\infty \sqrt{\alpha(k) \wedge H(t)  } \  dt < \infty  
 %\quad \text{where} \quad V(t)= \mathrm{Var}\left ( {\bf 1}_{X_1 \leq t}- {\bf 1}_{Y_1 \leq t}\right )
\end{equation}
(in fact, he proved the result  for a weaker version of the $\alpha$-mixing coefficient, the same as that used in \cite{DM17} for the central limit theorem). It follows from Section 5 of \cite{DM17}, that the condition \eqref{condalpha} is always less restrictive than \eqref{CunyCond}. 
\end{rem}

As a consequence of Theorem \ref{th1}, we get the compact law of the iterated logarithm. Let $K$ be the unit ball of the reproducing kernel Hilbert space (RKHS) associated with $\Gamma$, and $C$ be the image of $K$ by the ${\mathbb L}_1(dt)$ norm. The following corollary holds:

\begin{cor} \label{compact}
Assume that \eqref{condalpha} is satisfied. Then the sequence
$$
    \frac{\sqrt{n}}{\sqrt{2 \log \log n}}W_1(\mu_n, \mu)
$$
is almost surely relatively compact, with limit set $C$.
\end{cor}

The proof of Theorem \ref{th1} is based on two ingredients: a martingale approximation in ${\mathbb L}_1(dt)$, as in \cite{DM17}, and the following version of the bounded law of the iterated logarithm, which has an interest in itself.

\begin{prop} \label{LILW1}
Assume that \eqref{condalpha} holds, and let 
\begin{equation}\label{condalphabis}
V= \int_0^{\infty} \sqrt{ \sum_{k=0}^\infty  \left (\alpha(k) \wedge H(t)
\right ) } \ dt \, .
\end{equation}
Then, there exists a universal constant $\eta$ such that for any $\varepsilon >0$,
\beq \label{boundedLIL}
\sum_{n \geq 2} \frac{1}{n} {\mathbb P} \left ( \max_{1 \leq k \leq n} k W_1 (\mu_k , \mu)  > ( \eta V + \varepsilon ) \sqrt{n \log \log n }\right ) < \infty \, .
\eeq
\end{prop}

\begin{rem} {\it (The bivariate case)}.
Let $(X_i, Y_i)_{i \in {\mathbb Z}}$ be a stationary sequence of ${\mathbb R}^2$-valued random variables, and define the coefficients $\alpha(k)$ as in \eqref{defalpha}, with 
the two $\sigma$-algebras ${\mathcal F}_0= \sigma(X_i, Y_i, i \leq 0)$ and ${\mathcal G}_k= \sigma(X_i, Y_i, i \geq k)$.
Let $\mu_X$ (resp. $\mu_Y$) be the common distribution of the $X_i$'s (resp. the $Y_i$'s), and let
$$
\mu_{n,X} = \frac 1 n \sum_{k=1}^n \delta_{X_k} \quad \text{and} \quad 
\mu_{n,Y} = \frac 1 n \sum_{k=1}^n \delta_{Y_k}\, .
$$
Combining the arguments in \cite{BDM} and the proof of Theorem \ref{th1}, one can prove the following strong invariance principle for 
$
n \left (W_1(\mu_{n, X}, \mu_{n,Y})- W_1(\mu_X, \mu_Y) \right)
$.

%Our goal is to give sufficient conditions under which 
%$$
%n \left (W_1(\mu_{n, X}, \mu_{n,Y})- W_1(\mu_X, \mu_Y) \right)
%$$
%satisfies a strong invariance principle. To do so, we follow the approach in \cite{BDM}, based on the %decomposition (2.6) of that paper. Let $\varphi$ be the continuous function  from ${\mathbb L}_1(dt)$ %to ${\mathbb R}$ defined by
%\begin{equation}\label{phi}
 % \varphi(x)= \int \left (\text{sign}\{ F_X(t)-F_Y(t)  \}  \, x(t)  {\bf 1}_{F_X(t) \neq F_Y(t)}  +|x(t)|   
 %{\bf 1}_{F_X(t) = F_Y(t)}\right) dt  \, , 
%\end{equation}
%where $F_X$ (resp. $F_Y$) is the cumlative distribution function of $\mu_X$ (resp. $\mu_Y$). 
%Note that $\varphi$ is Lipschitz: $|\varphi (x) -\varphi(y)| \leq \|x-y\|_{{\mathbb L}_1}$. 
%Combining the arguments in \cite{BDM} and the proof of Theorem \ref{th1}, we easily obtain the %following result.
%\begin{thm}\label{th2}
Let $\varphi$ be the continuous function  from ${\mathbb L}_1(dt)$ to ${\mathbb R}$ defined by
\begin{equation*}\label{phi}
  \varphi(x)= \int \left (\text{sign}\{ F_X(t)-F_Y(t)  \}  \, x(t)  {\bf 1}_{F_X(t) \neq F_Y(t)}  +|x(t)|   
 {\bf 1}_{F_X(t) = F_Y(t)}\right) dt  \, , 
\end{equation*}
where $F_X$ (resp. $F_Y$) is the cumulative distribution function of $\mu_X$ (resp. $\mu_Y$). 
Assume that 
\begin{equation*}\label{alphaV2}
  \int_0^\infty  \sqrt{ \sum_{k=0}^\infty  \left ( \alpha(k) \wedge H_X(t)
\right ) } \ dt < \infty \quad \text{and}  \quad \int_0^\infty  \sqrt{ \sum_{k=0}^\infty  \left ( \alpha(k) \wedge H_Y(t)
\right ) } \ dt < \infty \, .
\end{equation*}
Then, enlarging the probability space if necessary, there exists a sequence of i.i.d.  ${\mathbb L}_1(dt)$-valued  centered Gaussian random variables $(Z_i)_{i \geq 1}$ with covariance function given by: for any $f, g \in {\mathbb L}_\infty (dt)$,
\begin{multline*}\label{covalpha2}
\widetilde \Gamma(f,g)=  \mathrm{Cov} \left(  \int f(t) Z_1(t) \ dt , \int g(t) Z_1(t) \ dt \right )\\
=\sum_{k \in {\mathbb Z}}
\iint f(t)g(s) \mathrm{Cov} (  {\bf 1}_{X_0 \leq t} -  {\bf 1}_{Y_0 \leq t},  {\bf 1}_{X_k \leq s} -  {\bf 1}_{Y_k \leq  s})  \ ds \, dt  \, ,
\end{multline*}
and such that 
$$
n \left (W_1(\mu_{n,X}, \mu_{n, Y}) - W_1 (\mu_X, \mu_Y) \right ) - \varphi \left ( \sum_{k=1}^n Z_k
\right ) = o(\sqrt{n \log \log n}) \quad
 \text{almost surely.}
$$ 
%is almost surely relatively compact, with limit set
%$\varphi(K)$ (where $\varphi$ is the function defined in \eqref{phi}). 
%\end{thm}

%As a consequence of Theorem \ref{th2}, we get the compact law of the iterated logarithm. Let $\tilde K$ be the unit ball of the RKHS associated with $\widetilde \Gamma$, and $\tilde C$ be the image of $\tilde K$ by $\varphi$. The following corollary holds:

%\begin{cor}\label{Lilbi}
%Assume that \eqref{alphaV2} holds. Then the sequence
%$$
%    \frac{\sqrt{n}}{\sqrt{2 \log \log n}}\left (W_1(\mu_{n,X}, \mu_{n, Y}) - W_1 (\mu_X, \mu_Y) \right )
%$$
%is almost surely relatively compact, with limit set $\tilde C$.
%\end{cor}
\end{rem}

\section{Rates of convergence of the empirical estimator of the Conditional Value at Risk}\label{ES}
\setcounter{equation}{0}
The Conditional Value at Risk at level $u \in (0,1]$  of a real-valued integrable random variable  $X$ ($CVaR_u(X)$) is a ``risk measure" (according to the definition of Acerbi and Tasche \cite{AT}), which is widely used in mathematical finance. It is sometimes called Expected Shortfall of Average Value at Risk. We refer to the paper \cite{AT} for a clear definition of that indicator, and for its relation with other well known measures, such as the Value at Risk, the Worst Conditional Expectation, the Tail Conditional Expectation... According to Acerbi and Tasche \cite{AT}, $CVaR_u(X)$ can be expressed as 
$$
   CVaR_u(X)= -\frac 1 u \int_0^u F_X^{-1} (x) dx  \, ,
$$
%The Expected Shortfall at level $u \in (0,1]$ of a real-valued and integrable random variable $X$ 
%($ES_u(X)$) has been introduced by Acerbi et al. \cite{ANS}, and is now widely used in mathematical %finance. It is a modification of the well known Conditional Value at Risk (CVaR$_u(X)$) of  $X$, which %satisfies some nice additional properties ($ES_u(X)$ is a ``Risk measure" according to the definition %in \cite{AT}). According to Acerbi and Tasche \cite{AT}, $ES_u(X)$ can be expressed as 
%$$
%   ES_u(X)= -\frac 1 u \int_0^u F_X^{-1} (x) dx  \, ,
%$$
where $F_X$ is the cumulative distribution function of the variable $X$, and $F_X^{-1}$ is its usual cadlag inverse: $F_X^{-1}(u)= \inf \{x \in  {\mathbb R} : F_X(x) \geq u \}$. 

Concerning the difference between the Conditional Value at Risk of two random variables $X$ and $Y$, the following elementary inequality  holds (see for instance \cite{Rio17}):
\begin{equation} \label{diff}
  \left |CVaR_u(X)-CVaR_u(Y)   \right | \leq \frac 1 u \int_0^1 | F_X^{-1}(x) -F_Y^{-1}(x)| dx = \frac 1 u W_1(\mu_X, \mu_Y) \, ,
\end{equation}
where $\mu_X$ (resp. $\mu_Y$) is the distribution of $X$ (resp. $Y$). 

Consider now the problem of estimating $CVaR_u(X)$ from the random variables $X_1, ...,  X_n$, where $(X_i)_{i \in {\mathbb Z}}$ is a stationary sequence of $\alpha$-mixing random variables with common distribution $\mu=\mu_X$. A natural estimator is then 
$$
  \widehat{CVaR}_{u,n} =  -\frac 1 u \int_0^u F_n^{-1} (x) dx \, ,
$$
where $F_n$ is the empirical distribution function based on $X_1, \ldots, X_n$. From \eqref{diff}, we get the upper bound 
$$
 \left |CVaR_u(X)-  \widehat{CVaR}_{u,n} \right | \leq \frac 1 u \int_0^1 |F_X^{-1}(x) -F_n^{-1}(x)| dx = \frac 1 u W_1(\mu_n, \mu) \, ,
$$
From Corollary \ref{compact}, we obtain the almost sure  rate of convergence of 
$ \widehat{CVaR}_{u,n}$: if  \eqref{condalpha} holds, then 
$$
     \limsup_{n \rightarrow \infty} \frac{\sqrt{n}}{\sqrt{2 \log \log n}}  \left |CVaR_u(X)-  \widehat{CVaR}_{u,n} \right | \leq \frac{\kappa(\Gamma)}{u } \ \  \text{almost surely},
$$
where $\kappa(\Gamma)$ is the largest value of the compact set $C$ of Corollary \ref{compact} (recall that the covariance function $\Gamma$ is defined in \eqref{covalpha}).  It is well known (see for instance Section 8 in \cite{LT}) that the constant $\kappa(\Gamma)$ can be expressed as
$$
   \kappa(\Gamma)= \sup_{f : \|f \|_\infty \leq 1} \left ({\mathrm{Var}} \left ( \int f(t) Z(t) dt  \ \right ) \right )^{1/2} \leq  \left \| \int |Z(t)| dt \right \|_2 \, ,
$$
where $Z$ is an ${\mathbb L}_1(dt)$-valued centered random variable with covariance function $\Gamma$. 
\section{Proofs}

\setcounter{equation}{0}

\subsection{Proof of Theorem \ref{th1}}
Let $(\Omega, {\mathcal A},  {\mathbb P})$ be the underlying probability space. By a standard argument, one may assume that $X_i=X_0 \circ T$, where $T: \Omega \mapsto \Omega$ is a bijective, bi-measurable transformation, preserving the probability ${\mathbb P}$. Let also ${\mathcal F}_i= \sigma(X_k, k \leq i)$. 

Let $Y_0(t)= {\mathbf 1}_{X_0 \leq t}-F(t)$, and $Y_k(t)=Y_0(t) \circ T^k= {\mathbf 1}_{X_k \leq t}-F(t)$. 
With these notations and the representation \eqref{Frep} one has that 
\begin{equation}\label{nW1sum}
n W_1(\mu_n, \mu)= \int \left | \sum_{k=1}^n Y_k(t) \right | \, dt \, .
\end{equation}
From Section 4 in \cite{DM17}, we know that, if \eqref{condalpha} holds, then
\begin{equation}\label{cobord}
  Y_0(t)= D_0(t) + A(t) - A(t) \circ T ,
\end{equation}
where $D_0$ is such that ${\mathbb E}(D_1(t)|{\mathcal F_{-1}})=0$ almost surely and $\int \|D_0(t)\|_2 \,  dt < \infty$, and $A$ is such that $\int \|A(t)\|_1 \,  dt < \infty$. Moreover, the covariance operator of $D_0$ is exactly $\Gamma$: for any $f, g \in {\mathbb L}_\infty (dt)$,
\begin{equation}\label{covD}
\Gamma(f,g)= \mathrm{Cov} \left(  \int f(t) D_0(t) \ dt , \int g(t) D_0(t) \ dt \right ) \, .
\end{equation}
Let $D_k(t)= D_0 \circ T^k$. From \eqref{cobord}, it follows that
\begin{equation}\label{main}
\sum_{k=1}^n Y_k= \sum_{k=1}^n D_k + A\circ T - A \circ T^n \, .
\end{equation}
From \cite[Proposition 3.3]{Cuny17}, we know that, enlarging the probability space if necessary, there exists a sequence of i.i.d. ${\mathbb L}_1(dt)$-valued centered Gaussian random variables $(Z_i)_{i \geq 1}$ with covariance function $\Gamma$ such that
\begin{equation}\label{StrongMart}
  \int \left |   \sum_{k=1}^n D_k(t) - \sum_{k=1}^n Z_k(t) \right | \, dt= o\left( \sqrt{n \log \log n}\right) \quad \text{almost surely}.
\end{equation}
Hence, the result will follow from \eqref{nW1sum}, \eqref{main} and \eqref{StrongMart} if we can prove that 
\begin{equation}\label{ascobord}
\lim_{n \rightarrow \infty} \frac{1}{\sqrt{n \log \log n}}\int \left | A(t)\circ T^n \right | \, dt =0 \quad \text{almost surely.}
\end{equation}

To prove \eqref{ascobord}, we start by considering the integral over $[-M,M]^c$, for $M>0$. Applying again  \cite[Proposition 3.3]{Cuny17}, we infer that 
\begin{equation}\label{Mc1}
\limsup_{n \rightarrow \infty} \frac{1}{\sqrt{2n \log \log n}} \int_{[-M,M]^c} \left | \sum_{k=1}^n D_k(t)\right | \, dt  \leq  \int_{[-M,M]^c} \|D_0(t)\|_2 \, dt \quad \text{almost surely.}
\end{equation}
Now, as will be clear from the proof, Proposition \ref{LILW1} also holds on the space ${\mathbb L}^1([-M,M]^c, dt)$, and implies that there exists a universal constant $\eta$ such that, for any positive $\varepsilon$, 
\begin{equation}\label{Mc2}
\limsup_{n \rightarrow \infty} \frac{1}{\sqrt{n \log \log n}} \int_{[-M,M]^c} \left | \sum_{k=1}^n Y_k(t)\right | \, dt  \leq \varepsilon +  \eta  \int_M^\infty   \sqrt{ \sum_{k=0}^\infty \min \left\{\alpha(k), H(t)
\right \} } \ dt  \quad \text{almost surely.}
\end{equation}
From \eqref{Mc1} and \eqref{Mc2}, we infer that
$$
\lim_{M\rightarrow \infty} \limsup_{n \rightarrow \infty} \frac{1}{\sqrt{n \log \log n}} \int_{[-M,M]^c} \left | A(t)\circ T^n \right| \, dt =0 \quad \text{almost surely.}
$$
Hence the proof of \eqref{ascobord} will be complete if we prove that, for any $M>0$,
\begin{equation}\label{laststep}
\limsup_{n \rightarrow \infty} \frac{1}{\sqrt{n \log \log n}} \int_{-M}^M \left | A(t)\circ T^n \right| \, dt =0 \quad \text{almost surely.}
\end{equation}
To prove \eqref{laststep}, we work in the space ${\mathbb H}={\mathbb L_2}([-M, M], dt)$, and we denote by $\|\cdot \|_{\mathbb H}$ and $\langle \cdot, \cdot \rangle$ the usual norm and scalar product on ${\mathbb H}$. Since ${\mathbb E}(\|D_0\|^2_{\mathbb H}) < \infty$, we know from \cite{Cuny17}
 that $\sum_{k=1}^n D_k $ satisfies the compact law of the iterated logarithm in ${\mathbb H}$. Since $\sum_{k \geq 0} \alpha(k) < \infty$  and $Y_0$ is bounded in ${\mathbb H}$, we infer from 
 \cite{DM10}  that  $\sum_{k=1}^n Y_k $ satisfies also the compact law of the iterated logarithm in 
${\mathbb H}$.

Now, arguing exactly as in the end of the proof of  \cite[Theorem 4]{DM10}, one has: for any $f$ in ${\mathbb H}$
\begin{equation}\label{DM10}
\lim_{n \rightarrow \infty} \frac{\langle f, A \circ T^n \rangle}{\sqrt{n \log \log n} }= 0 \quad \text{almost surely.}
\end{equation}
Let $(e_i)_{i \geq 1}$ be a complete orthonormal basis of ${\mathbb H}$ and 
$P_N(f)=\sum_{k=1}^N \langle f, e_k \rangle e_k$ be the projection of $f$ on the space spanned by the first $N$ elements of the basis. From \eqref{DM10}, we get that 
\begin{equation}\label{PN}
\lim_{n \rightarrow \infty} \frac{P_N(A \circ T^n) }{\sqrt{n \log \log n} }= 0 \quad \text{almost surely.}
\end{equation}
On another hand, applying again  \cite[Proposition 3.3]{Cuny17} (as done in \eqref{Mc1}), we get 
\begin{equation}\label{Mart}
\lim_{N \rightarrow \infty} \limsup_{n \rightarrow \infty}  \frac{1}{\sqrt{n \log \log n}}  \left \|  (I-P_N) \left( \sum_{k=1}^n D_k  \right )\right \|_{\mathbb H}= 0 \quad \text{almost surely,}
\end{equation}
and applying \cite[Theorem 4]{DM10}, 
\begin{equation}\label{Y}
\lim_{N \rightarrow \infty} \limsup_{n \rightarrow \infty} \frac{1}{\sqrt{n \log \log n}} \left \|  (I-P_N) \left( \sum_{k=1}^n Y_k  \right )\right \|_{\mathbb H}= 0 \quad \text{almost surely.}
\end{equation}
From \eqref{main}, \eqref{Mart} and \eqref{Y}, we infer that 
\begin{equation*}
\lim_{N \rightarrow \infty} \limsup_{n \rightarrow \infty} \frac{\|(I-P_N)A \circ T^n \|_{\mathbb H}}{\sqrt{n \log \log n} }= 0 \quad \text{almost surely, }
\end{equation*}
which, together with  \eqref{PN}, implies \eqref{laststep}. The proof of Theorem \ref{th1} is complete. $\diamond$

\subsection {Proof of Proposition \ref{LILW1}}
For any $n \in {\mathbb N}$, let us introduce the following notations:
  \begin{equation*}
  R(u)=\min\{q\in {\mathbb N}^* \tq \alpha(q)\leq u\}  Q(u)
\quad   \text{and} \quad 
R^{-1}(x)=\inf\{u\in [0,1] \tq R(u)\leq x\}\,.
  \end{equation*}
For  a positive real $a$ that will be specified later,   let
\beq \label{defvM}
m_n = a \sqrt{\frac{n}{\log \log  n }}\, , \quad  v_n = R^{-1}(m_n) \, , \quad  M_n=Q(v_n) \, .
\eeq
For any $M>0$, let  $ g_M(y) = (y \wedge M) \vee (-M) $.  For any integer $i$, define
\beq \label{detrun}
X_i' = g_{M_n}(X_i )\,  \mbox{ and } \, X_i'' = X_i -X_i' \, .
\eeq
We first recall that, by the dual expression of $W_1 (\mu_n , \mu)$, 
$$
n W_1 (\mu_n , \mu)   = \sup_{f \in \Lambda_1} \sum_{i=1}^n \left ( f(X_i) - \E(f(X_i)) \right ) \, .
$$
where $\Lambda_1$ is the set of Lipschitz functions such that $|f(x)-f(y)| \leq |x-y|$. 
Hence,
$$
n W_1 (\mu_n , \mu)  
\leq
\sup_{f \in \Lambda_1} \sum_{i=1}^n
\left ( f(X'_i) - \E(f(X'_i)) \right )
+
\sup_{f \in \Lambda_1} \sum_{i=1}^n
\left ( f(X_i) - f(X_i') -  \E(f(X_i)-f(X_i')) \right )  \, .
$$
Therefore, setting,
\[
F_n' (t) = \frac{1}{n} \sum_{k=1}^n {\bf 1}_{\{X'_k \leq t \} } \quad  \text{and} \quad   F' (t) = {\mathbb P} ( X_1' \leq t) \, ,
\]
and noticing that
\[
k  \Vert F_k' - F' \Vert_1 = \sup_{f \in \Lambda_1} \sum_{i=1}^k
  \big ( f(X_i') -  \E(f(X_i') \big )  \, ,
\]
we get
\beq \label{decW1Wprime}
 \max_{1 \leq k \leq n} k W_1 (\mu_k , \mu)    \leq  \max_{1 \leq k \leq n} k  \Vert F_k' - F' \Vert_1
+ \sum_{i=1}^n ( |X_i''| + \E( |X_i''|) \, .
\eeq
Now, note that
\begin{multline*}
\sum_{n \geq 2 } \frac{1}{ \sqrt {n \log \log n}  } \E( |X_n''|)  \leq
\sum_{n \geq 2 } \frac{1}{ \sqrt {n \log \log n}  } \int_{0}^{+ \infty} {\mathbb P } \big ( |X_0| {\bf 1}_{|X_0| >Q (v_n) }  > t  \big )  dt \\
 \leq
\sum_{n \geq 2 } \frac{1}{ \sqrt {n \log \log n}  } \int_{Q(v_n)}^{+ \infty} H(t)   dt 
 \leq
\sum_{n \geq 2 } \frac{1}{ \sqrt {n \log \log n}  } \int_{0}^{v_n} Q(u) du 
\\ \leq  \sum_{n \geq 2 } \frac{1}{ \sqrt {n \log \log n}  } \int_{0}^{1} Q(u) {\mathbf 1}_{m_n \leq R(u)} du  \ll \int_0^1 R(u) Q(u) du  \, .
\end{multline*}
But, according to Propositions 5.1 and 5.2 in \cite{DM17}, condition \eqref{condalpha} implies that
\begin{equation}\label{condalpha2}
\int_0^1 R (u) Q(u) du < \infty \, .
\end{equation}
Hence, to prove \eqref{boundedLIL} it suffices to show that there exists an universal constant $\eta$ such that for any $\varepsilon >0$,
\beq \label{boundedLIL2}
\sum_{n \geq 2} \frac{1}{n} {\mathbb P} \left (  \max_{1 \leq k \leq n} k  \Vert F_k' - F' \Vert_1  > \eta V \sqrt{n \log \log n }\right ) < \infty \, .
\eeq
For this purpose, let 
\beq \label{defofq}
q_n = \min\{k\in {\mathbb N}^* \tq \alpha(k)\leq v_n\} \wedge n \, .
\eeq
Since $R$ is right
continuous, we have $R(R^{-1}(w))\leq w$ for any $w$, hence
\beq \label{restq1}
q_nM_n = R(v_n) = R (R^{-1}(m_n)) \leq m_n \, .
\eeq
Assume first that $q_n=n$. Bounding  $f(X'_i) - \E(f(X'_i)) $ by $2M_n$, we obtain
\begin{equation}\label{vB1}
 \max_{1 \leq k \leq n} k  \Vert F_k' - F' \Vert_1  \leq 2 n M_n  = 2q_nM_n  \leq 2 m_n \, .
\end{equation}
Taking into account the definition of $m_n$, it follows that there exists  $n_0$ depending on $a$, $V$ and $\eta$, such that for any $n \geq n_0$, $8m_n \leq  \kappa V  \sqrt{ n \log \log n} $.  This proves the proposition  in the case where $q_n=n$.

From now on, we assume that $q_n<n$. Therefore $q_n =  \min\{k\in {\mathbb N}^* \tq
\alpha (k)\leq v_n\} $ and then $\alpha(q_n)\leq v_n$. For any integer $i$, define
\[ 
U_i(t) =  \sum_{k=(i-1)q_n+1}^{iq_n} \left(\I_{X_k' \leq t} - {\mathbb E}\left  (\I_{X_k' \leq t} \right ) \right ) \, .
\]
and notice that
\[
 \max_{1 \leq k \leq n} k  \Vert F_k' - F' \Vert_1  \leq  2q_nM_n  + \int_{-M_n}^{M_n} \max_{1 \leq j \leq [n/q_n]} \left |  \sum_{i=1}^j U_{i} (t) \right | dt  \, .
\]

Let $k_n = [n /q_n]$.  For any $t$, applying  Rio's coupling lemma (see
 \cite[Lemma 5.2]{Rio00}) recursively,  we can construct random variables $(U_i^*(t) )_{1 \leq i \leq k_n }$ such that 
 \begin{itemize}
  \item $U_i^*(t) $ has the same distribution as $U'_i$ for all $1 \leq i \leq k_n $,\\
  \item the random variables $(U_{2i}^*(t) )_{2 \leq 2i \leq k_n }$ are independent, 
     as well as the random variables $(U_{2i-1}^*(t) )_{1 \leq 2i-1 \leq k_n }$,\\
  \item we can suitably control $\|U_i (t) - U_i^*(t) \|_1$ as follows: for any $i \geq 1$,
\beq \label{couplingineRio}
  \|U_i (t) - U_i^*(t) \|_1 \leq 4 q_n \alpha (q_n) \, .
\eeq
\end{itemize}
Substituting $U^*_i(t)$ to $U_i(t)$, we
obtain
\begin{multline} \label{decavecmart}
 \max_{1 \leq k \leq n} k  \Vert F_k' - F' \Vert_1  \leq  2q_nM_n
+ \max_{2 \leq 2j \leq
 [n/q_n]
 } \left \vert \sum_{i=1}^j U^*_{2i} (t) \right \vert \! \\ + \!\max_{1 \leq 2j-1 \leq
 [n/q_n]
 } \left \vert \sum_{i=1}^j U^*_{2i -1} (t) \right \vert
 \!+\!
 \sum_{i=1}^{ [n/q_n]
 } \vert  U_{i} (t) -U^*_{i} (t) \vert \, .
\end{multline}
Therefore, setting $\kappa=\eta/4$, for $n \geq n_0$,
\begin{equation} \label{decavecmartproba}
{\mathbb P} \left (  \max_{1 \leq k \leq n} k  \Vert F_k' - F' \Vert_1  \geq 4V \kappa  \sqrt{ n \log \log n} \right )
\leq  I_1(n) + I_2(n) + I_3(n) \, ,
\end{equation}
where
\begin{align*}
I_1(n)&={\mathbb P} \left (  \int_{-M_n}^{M_n}  \sum_{i=1}^{ [n/q_n]
 } \vert  U_{i} (t) -U^*_{i} (t) \vert  \ dt \geq  V \kappa \sqrt{ n \log \log n} \right )  \\
I_2(n) & = {\mathbb P} \left ( \int_{-M_n}^{M_n} \max_{2 \leq 2j \leq
 [n/q_n]
 } \left \vert \sum_{i=1}^j U^*_{2i} (t) \right \vert \ dt  \geq  V \kappa \sqrt{ n \log \log n} \right )  \\
I_3(n) & =   {\mathbb P} \left (  \int_{-M_n}^{M_n} \max_{1 \leq 2j-1 \leq
 [n/q_n]
 } \left \vert \sum_{i=1}^j U^*_{2i-1} (t)\right \vert \ dt  \geq  V \kappa \sqrt{ n \log \log n} \right ) \, .
\end{align*}
Using Markov's inequality and \eqref{couplingineRio}, we get
\begin{equation*} I_1(n) \ll  \frac{n}{\sqrt{ n \log \log n} } M_n \alpha (q_n )   \ll  \frac{n}{\sqrt{ n \log \log n} } v_n Q(v_n) \ll \frac{n}{\sqrt{ n \log \log n} } \int_0^{R^{-1} (m_n) } Q(u) du    \, .
\end{equation*}
Hence, by \eqref{condalpha2}, 
\[
\sum_{n \geq 2}  \frac{1}{n} I_1 (n) \ll  \sum_{n \geq 2} \frac{1}{\sqrt{ n \log \log n} } \int_0^{R^{-1} (m_n) } Q(u) du  \ll \int_0^1 R(u) Q(u) du  < \infty \, .
\]

To handle now the term $I_2(n)$ (as well as   $I_3(n)$)  in the decomposition \eqref{decavecmartproba}, we shall use again  Markov's inequality but this time at the order $p \geq 2$. Hence for $p \geq 2$, taking into account the stationarity, we get
\[
I_2(n)  \leq \frac 1 { ( V \kappa)^p (n \log \log n )^{p/2} } \left ( \int_{-Q(v_n)}^{Q(v_n)}
\left \| \max_{2 \leq 2j \leq
 [n/q_n] } \left | \sum_{i=1}^j \tilde U_{2i} (t) \right |
\right  \|_p dt   \right )^p  \, .
\]
Applying Rosenthal's inequality (see for instance \cite[Theorem 4.1]{Pi94}) and taking into account the stationarity, there exist two  positive universal constants $c_1$ and $c_2$ not depending on $p$ such that
\beq \label{RosenthalU}
\left \| \max_{2 \leq 2j \leq
 [n/q_n] } \left | \sum_{i=1}^j U^*_{2i} (t) \right |
\right  \|_p^p \leq c_1^p p^{p/2}  (n/q_n)^{p/2} \Vert  U_2 (t)  \Vert^p_2 + c_2^p p^p (n/q_n) \Vert  U_2 (t)  \Vert^p_p := J_1(t) + J_2 (t)  \, .
\eeq
Using similar arguments as to handle the quantity $I_2(n)$ in the proof of \cite[Proposition 3.4]{DM17}, we have
%\marginpar{donner des explications}
\begin{multline}\label{intQvn}
 \int_{-Q(v_n)}^{Q(v_n)}     \Vert  U_2 (t)  \Vert_2  dt =   \int_{-Q(v_n)}^{Q(v_n)}   \left (  {\rm Var}  \left( \sum_{i=1}^{q_n} {\bf 1}_{\{ X_i' \leq t \} }\right )   \right )^{1/2}  dt \\
 \leq 2 \sqrt{2} \sqrt{q_n}  \int_{0}^{Q(v_n)}    \left (   \sum_{k=0}^{q_n-1} \alpha(k) \wedge   H(t)  \right )^{1/2} dt \leq 2 V   \sqrt{2 q_n} \, .
\end{multline}
Hence
\[
\sum_{n \geq 2}  \frac{1}{n ( V \kappa)^p (n \log \log n )^{p/2} } \left (   \int_{-Q(v_n)}^{Q(v_n)}     J_1(t)^{1/p}  dt \right )^p  \leq  \sum_{n \geq 2}  \frac{ (2 \sqrt{2} c_1 \sqrt{p} )^p}{n  \kappa^p ( \log \log n )^{p/2} }   \, .
\]
Let now
\[
p=p_n =  \max \{ c \log \log n,   2 \},
\]
where $c$ will be specified later.
Set $n_1 = \min \{ n \geq 2 \, : \, c \log \log n \geq 2 \} $.  It follows that
\[
\sum_{n \geq n_1}  \frac{1}{n ( V \kappa)^p (n \log \log n )^{p/2} } \left (   \int_{-Q(v_n)}^{Q(v_n)}     J_1(t)^{1/p}  dt \right )^p  \leq  \sum_{n \geq n_1}  \frac{1}{n}  \left (   \frac{ 2 c_1 \sqrt{2c} }{\kappa} \right ) ^{c \log \log n}   \, ,
\]
which is finite provided we take $\kappa$ such that $ \frac{ 2 c_1 \sqrt{2c} }{\kappa} = \alpha^{-1} $ with $\alpha >1$ and $c > (\log \alpha)^{-1}$.

On another hand, proceeding as in \eqref{intQvn}, we deduce that, for any $t >0$,
\begin{multline*}
 \Vert  U_2 (t)  \Vert_p^p = \left \Vert  \sum_{i=1}^{q_n} \left (  {\bf 1}_{\{ X_i' \leq t \} } - {\mathbb P} ( X_i' \leq t ) \right ) \right \Vert_p^p
 \leq q_n^{p-2}  \left \Vert  \sum_{i=1}^{q_n} \left (  {\bf 1}_{\{ X_i' \leq t \} } - {\mathbb P} ( X_i' \leq t ) \right ) \right \Vert_2^2  \\
 \leq  2 q_n^{p-1}   \sum_{k=0}^{q_n-1} ( \alpha(k) \wedge   H(t)  ) \, .
\end{multline*}
In addition
\begin{multline*}
\int_0^{Q(v_n)}  \left (  \sum_{k=0}^{q_n-1} \alpha(k) \wedge   H(t)  \right )^{1/p} dt =  \int_0^{Q(v_n)}  \left (  \int_0^{H(t) }  (\alpha^{-1} (u) \wedge q_n )du    \right )^{1/p} dt  \\
\leq \int_0^{Q(v_n)}  \left (  v_nq_n  +  \int_{v_n}^{H(t) }  (\alpha^{-1} (u) \wedge q_n ) du    \right )^{1/p} dt  \, . \end{multline*}
Note that $u < H(t) \iff t < Q(u)$.  Consequently $u < H(t) $ implies that $Q^{-2} (u) < t^{-2}$. Hence
\begin{multline*}
\int_0^{Q(v_n)}  \left (  \sum_{k=0}^{q_n-1} \alpha(k) \wedge   H(t)  \right )^{1/p} dt
\\ \leq (v_n q_n)^{1/p}  Q(v_n)   + \int_0^{Q(v_n)} \left ( t^{-2}  \int_{v_n}^{H(t) }   (\alpha^{-1} (u) \wedge q_n ) Q^2(u) du    \right )^{1/p}  \\
\leq (v_n q_n)^{1/p}  Q(v_n) +  \left (  \int_{v_n}^{1}  R(u) Q (u) du    \right )^{1/p}   \int_0^{Q(v_n)} t^{-2/p} dt  \\
\leq  (v_n q_n)^{1/p}  Q(v_n) + \left (  \int_{0}^{1}  R(u) Q (u) du    \right )^{1/p} p(p-2)^{-1} Q(v_n)^{1-2/p}  \, .
 \end{multline*}
Set $n_2= \min \{ n \geq 2 \, : \, c \log \log  n \geq 4 \} $.  It follows that
\begin{multline*}
\sum_{n \geq n_2}  \frac{1}{n ( V \kappa)^p (n \log \log n )^{p/2} } \left (   \int_{-Q(v_n)}^{Q(v_n)}     J_2(t)^{1/p}  dt \right )^p \\ \leq  2  \sum_{n \geq n_2}  \frac{ ( 4 c_2  p )^p}{ (\kappa V)^p (n \log \log n )^{p/2} }  q_n^{p-2}
\left \{ v_n q_n Q^p(v_n)  +  2^p Q(v_n)^{p-2}    \int_{0}^{1}  R(u) Q (u) du  \right \}   \, .
 \end{multline*}
Note that
\[
v_n q_n Q^2(v_n)   = v_n \alpha^{-1} (v_n) Q^2(v_n)  \leq  \int_{0}^{1}  R(u) Q (u) du \, .
\]
Hence, since $q_nM_n  \leq m_n$, we get
\begin{multline*}
\sum_{n \geq n_2}  \frac{1}{n ( V \kappa)^p (n \log \log n )^{p/2} } \left (   \int_{-Q(v_n)}^{Q(v_n)}     J_2(t)^{1/p}  dt \right )^p \\ \leq  4  \int_{0}^{1}  R(u) Q (u) du \sum_{n \geq n_2}  \frac{ ( 8 c_2  p )^p}{ (\kappa V)^p (n \log \log n )^{p/2} }  m_n^{p-2} \\
\leq  4 a^{-2}  \int_{0}^{1}  R(u) Q (u) du \sum_{n \geq n_2}   \left ( \frac{ 8 a c_2   c  }{ \kappa V }  \right )^p \frac{\log \log n }{n } \, ,
 \end{multline*}
which is finite by taking into account \eqref{condalpha2}, and if we choose $a =  (  c_1 \kappa V)/(2c_2 \sqrt{2c})$. Indeed, in this case,
\[
 \frac{ 8 a c_2   c  }{ \kappa V } =   \frac{ 2 c_1 \sqrt{2c} }{\kappa}  \times   \frac{  2a c_2 \sqrt{2c} }{ c_1 \kappa V }= \alpha^{-1} \, .
\]
This ends the proof of the proposition. $\diamond$

\end{document}